\documentclass{acm_proc_article-sp}
\usepackage{amssymb}
\usepackage{latexsym}
\newtheorem{theorem}{Theorem}
\newtheorem{prop}{Proposition}
\newtheorem{lemma}{Lemma}
\def\leqs{{\leqslant}}
\def\geqs{{\geqslant}}

\def\Z{{\bf Z}}
\def\R{{\bf R}}
\def\C{{\bf C}}
\def\Pr{{\bf P}}
\def\CC{{\cal C}}
\def\Re{\mathop{\rm Re}\nolimits}
\def\PGL{{\rm PGL}}
\def\Pic{{\rm Pic}}
\def\vph{\varphi\0}
\def\ra{\rightarrow}

\def\0{^{\phantom0}}
\def\9{_{\phantom9}}

\def\GV{Gilbert-\mbox{\kern-.16em}Varshamov}
\def\DV{Drinfeld-\mbox{\kern-.16em}Vl\u{a}du\c{t}}
\def\zetaC{\zeta\0_C}

\begin{document}
\conferenceinfo{STOC'01,} {July 6-8, 2001, Hersonissos, Crete, Greece.}
\CopyrightYear{2001}
\crdata{1-58113-349-9/01/0007}

\title{Excellent nonlinear codes from modular curves}

\numberofauthors{1}
\author{
\alignauthor Noam D. Elkies\\
       \affaddr{Department of Mathematics}\\
       \affaddr{Harvard University}\\
       \affaddr{Cambridge, MA 02138 USA}\\
       \email{elkies@math.harvard.edu}
}
\date{30 October 2000}
\maketitle
\begin{abstract}
We introduce a new construction of error-correcting codes from
algebraic curves over finite fields.  Modular curves of genus
$g\ra\infty$ over a field of size $q_0^2$ yield nonlinear codes more
efficient than the linear Goppa codes obtained from the same curves.
These new codes now have the highest asymptotic transmission rates
known for certain ranges of alphabet size and error rate.
Both the theory and possible practical use of these new record codes
require the development of new tools.  On the theoretical side,
establishing the transmission rate depends on an error estimate for
a theorem of Schanuel applied to the function field of an
asymptotically optimal curve.  On the computational side,
actual use of the codes will hinge on the solution of new problems
in the computational algebraic geometry of curves.

\end{abstract}


\section{Prologue}
In this section we first review the construction and properties
of Goppa codes, to put our work in its context.
We then define our new nonlinear codes
and give lower bounds on their minimal distance.
We conclude this section
by stating lower bounds on the size of our codes
and comparing our codes' parameters with those of Goppa codes.
In the next section we prove the bounds claimed in the Introduction.
In the final section we discuss theoretical and computational questions
raised by our construction, and show how to solve these problems
for the nonlinear codes obtained from rational curves.

\subsection{Review: algebro-geometric (Goppa) codes}

Fix a finite field~$k$ of $q=p^\alpha$ elements.
Let $C$\/ be a projective, smooth, irreducible
algebraic curve of genus~$g$ defined over~$k$,
with $N$\/ rational points.  
To any divisor~$D$\/ on~$C$\/ of degree $<N$,
Goppa (\cite{Goppa}, see also \cite{TV})
regards the space of sections of~$D$\/
as a linear $[N,r,d]$ code with alphabet~$k$,
for some $d \;\geqs\; N-\deg(D)$
(because a nonzero section of~$D$\/ has at most $\deg(D)$ zeros)
and $r \geqs \deg(D)-g+1$ (by the Riemann-Roch theorem).
Thus the transmission rate $R=r/N$\/
and the error-detection rate $\delta=d/N$\/ of Goppa's codes
are related by
\begin{equation}
\label{eq:Goppa}
R + \delta > 1 - \frac{g}{N} .
\end{equation}
This lower bound improves as $N/g$ increases.
How large can $N/g$ get as $g\ra\infty$? An upper bound is
\begin{equation}
\label{eq:DV}
N < (q^{1/2}-1+o(1)) g
\end{equation}
(\DV~\cite{DV}).
We say a curve of genus $g\ra\infty$ is ``asymptotically optimal''
if it has at least $(q^{1/2}-1-o(1)) \, g$ rational points over~$k$.
If $\alpha$ is even, i.e., if $q_0\0 := \sqrt q$ is an integer,
then modular curves of various flavors ---
classical (elliptic), Shimura, or Drinfeld --- attain
\begin{equation}
\label{eq:modular}
N \ \geqs\ (q_0\0-1)(g-1) = (q^{1/2}-1-o(1)) \, g
\end{equation}
\cite{Ihara, TVZ}, and are thus asymptotically optimal.
Therefore if $q=q_0^2$
there exist arbitrarily long linear codes over~$k$\/ with
\begin{equation}
\label{eq:Goppamax}
R + \delta > 1 - \frac{1}{q_0\0-1} - o(1),
\end{equation}
and this is the best that can be obtained from~(\ref{eq:Goppa}).
Once \hbox{$q_0\0\;\geqs\;7$,} these codes improve on the \GV\ bound
for suitable $R,\delta$.

Actual construction of these codes requires explicit equations for~$C$.
The definitions of modular curves do not readily yield useful equations,
but in recent years many families of modular curves have been given
by $O(\log g)$ explicit equations in $O(\log g)$ variables,
each equation of degree $O(\log g)$.
See \cite{E:tower} for classical and Shimura curves,
\cite{E:shimura} for further Shimura curves, and
\cite{GS:drinfeld, GS:more_drinfeld, E:drinfeld, TV}
for Drinfeld modular curves.\footnote{
  Actually the equations in \cite{TV} are in two variables but
  of degree exponential in~$\log g$; but they are easily put in
  an equivalent form of degree $O(\log g)$ by introducing
  $O(\log g)$ more variables.
  }
Using the resulting codes for error-resistant communication
also requires polynomial-time decoding
of any word at distance $<d/2$ from a codeword;
this and more has also been recently accomplished~\cite{GuSu,ShWa}.

\subsection{The new nonlinear codes}

The Goppa codes generalize the Reed-Solomon codes,
which are the special case where $C$\/ is a projective line
$\Pr^1$ (so $g=0$).  In this special case, the Goppa code
can be identified with the space of polynomials of degree
at most~$\deg(D)$ in one variable, interpreted as words
by evaluation at each element of~$k$.\footnote{
  More precisely, the Goppa codes for $g=0$ are extended
  Reed-Solomon codes of length $q+1$, with one coordinate
  for each element of~$k$, and an additional coordinate
  for the leading coefficient, corresponding to
  evaluation at the point at infinity of~$\Pr^1$.
  }
Our new idea is to replace these polynomials by rational functions
of bounded degree, say degree $\leqs\; h$.  Since a rational function
of degree~$\leqs\; h$\/ is determined by two polynomials of degree
$\leqs\; h$, we expect that $h$\/ will play a role comparable to
{\em half}\/ the degree of the divisor~$D$\/ used to construct
a Goppa code.  The notions of a rational function and its degree
extend to curves~$C$\/ of arbitrary genus.  Given $C$\/ with
$N$\/ rational points, we thus define $\CC_0(h)$ for any $h<N/2$
as follows: $\CC_0(h)$ consists of the rational functions~$f$\/ on~$C$,
defined over~$k$, such that $\deg(f) \;\leqs\; h$.
To give $\CC_0(h)$ the structure of an error-correcting code,
choose an enumeration $(P_1,\ldots,P_N)$
of the $k$-rational points of~$C$,
and identify $f$\/ with the $N$\/-tuple
\begin{equation}
\label{eq:Ntuple}
\bigl( f(P_1),f(P_2),\ldots,f(P_N) \bigr)
\end{equation}
of values of~$f$\/ at points of~$C$.  Since $f$\/ may have poles
on some $P_i$, some $f(P_i)$ values may be $\infty$.
Thus the alphabet for our new code $\CC_0(h)$ is not a finite field
but a set of size $q+1$, the projective line
$k\cup\{\infty\}=\Pr^1(k)$ over the finite field~$k$.
In other words, we are identifying a function~$f$\/
with its graph as a map from~$C$\/ to~$\Pr^1(k)$,
just as a polynomial in the Reed-Solomon code
was identified with its graph as a map from~$k$ to~$k$.
It is readily seen (Prop.~\ref{prop:distance} below) that if $f_1,f_2$
are distinct rational functions of degrees $h_1,h_2$ on~$C$\/
then $f_1(P)=f_2(P)$ holds for at most $h_1+h_2$ points~$P$\/ of~$C$.
Therefore $\CC_0(h)$ has minimal distance at least $N-2h$.
In particular, since we assume $h<N/2$,
different functions of degree $\leqs\; h$\/
yield different words in $\CC_0(h)$.

More generally, let $D$\/ be a divisor {\em of degree zero} on~$C$.
For each $h<N/2$ we define $\CC_D(h)$ to be the set of rational sections
of degree $\leqs\; h$\/ of the line bundle~$L_D$ associated to~$D$.
That is, $\CC_D(h)$ consists of the zero function
together with the nonzero rational functions~$f$\/ on~$C$\/
whose divisor $(f)$ is of the form $E-D$\/ for some divisor~$E$\/
whose positive and negative parts each have degree at most~$h$.
To give~$\CC_D(h)$ the structure of an error-correcting code,
choose for each $k$-rational point $P_i$ of~$C$\/ a rational function
$\vph_i$ whose divisor has the same order at~$P_i$ as~$D$, and
identify each $\CC_D(h)$ with the $N$-tuple
\begin{equation}
\label{eq:NtupleD}
\bigl( (\vph_1 f)(P_1),(\vph_2 f)(P_2),\ldots,(\vph_N f)(P_N) \bigr)
\in (\Pr^1(k))^N .
\end{equation}
Different choices of $\vph_i$ yield isomorphic codes
(Lemma~\ref{lem:indep} below).
In particular, if $D=0$ we recover our earlier definition of $\CC_0(h)$
by setting each $\vph_i=1$.  We shall see in Prop.~\ref{prop:distance}
that here, too, any two distinct rational sections of degrees $h_1,h_2$
agree on at most $h_1+h_2$ points,
so $\CC_0(h)$ has minimal distance at least $N-2h$,
and $f$\/ can be recovered uniquely
from the $N$-tuple (\ref{eq:NtupleD}).
Linearly equivalent divisors yield isomorphic codes
(Lemma~\ref{lem:lineq}),
so $D$\/ can be regarded as a degree-zero divisor
modulo linear equivalence,
i.e., as an element of the {\em Jacobian} $J_C$ of~$C$.

\subsection{Size of the codes; comparison with Goppa}

Let $M(h,C)$ be the average size of $\CC_D(h)$
as $D$\/ varies over~$J_C$:
\begin{equation}
\label{eq:Mdef}
M(h,C) := \frac1{\#(J_C)} \sum_{D\in J_C} \#(\CC_D(h)).
\end{equation}
We shall show (Thm.~\ref{thm:rho0})
that if $C$\/ is an asymptotically optimal curve then,
for each
\begin{equation}
\label{eq:rho0}
\rho > \frac{2q}{q^2-1},
\end{equation}
the estimate
\begin{equation}
\label{eq:Masymp}
M(h,C) = \left(\frac{q+1}{q}\right)^{N+o_\rho(N)} q^{2h-g}
\end{equation}
holds as long as $2h/N > \rho$.  The threshold~(\ref{eq:rho0})
is low enough to allow all ratios $h/N$\/
for which the estimate~(\ref{eq:Masymp}) exceeds~$1$.
In particular, if $2h \;\geqs\; g$, our codes have on average
\begin{equation}
\bigl( (q+1)/q \bigr)^{N+o(N)}
\label{Mratio}
\end{equation}
times as many words as the Goppa codes
of the same length and designed minimal distance
must have by Riemann-Roch.
With a somewhat longer argument we show (Thm.~\ref{thm:rho1})
that the same estimate holds for each individual $\CC_D(h)$,
but with a higher threshold $\rho\0_1(q)$
defined below (equations \ref{eq:B1},\ref{eq:rho_1}).

We cannot simply conclude that our codes transmit asymptotically
$\log \bigl( (q+1)/q \bigr)$ more bits per letter than Goppa's, because
our alphabet size is larger by~$1$ than that of the Goppa codes.
A direct comparison would require
Goppa codes over a field of $q+1$ elements.
But it is rare that $q$ and $q+1$ are both prime powers
(one of them must be a power of~$2$,
the other a Mersenne or Fermat prime);
and they can never both be squares.
Nevertheless we claim that a fair comparison can be made,
and shows our codes to be better in a range of parameters
that includes all the Goppa codes that improve on~\GV.

We base this claim on two observations.
First, if a code over an alphabet of $q+1$ letters
is as good as a Goppa code, its parameters should obey the relation
obtained by extrapolating (\ref{eq:Goppamax}) to an alphabet
of size $q+1$, that is,
\begin{equation}
\label{eq:Goppaq+1}
R + \delta > 1 - \frac1{\sqrt{q+1}-1} - o(1).
\end{equation}
By (\ref{eq:Masymp}), our codes' parameters satisfy
\begin{equation}
\label{eq:Newmax}
\frac{\log(q+1)}{\log q} R + \delta
> 1 - \frac{1}{q_0\0-1} + \frac{\log \frac{q+1}{q}} {\log q} - o(1).
\end{equation}
This improves on (\ref{eq:Goppaq+1}) as long as
\begin{equation}
\label{eq:1-R}
1-R > \frac
  { \frac{1}{\sqrt{q}-1} - \frac{1}{\sqrt{q+1}-1} }
  { \log(\frac{q+1}{q}) / \log q }
= \frac{ \frac12 q^{-1/2} + O(q^{-1}) } {\log q}.
\end{equation}
This condition holds for all $(R,\delta)$ for which
(\ref{eq:Goppaq+1}) is better than the \GV\ bound.

For a second approach,
instead of extrapolating Goppa codes to alphabets of size $q+1$,
we degrade our codes by artificially reducing the alphabet size to~$q$.
To do this, we choose for each $i=1,\ldots,N$ a forbidden letter
$a_i \in \Pr^1(k)$, and consider only words $w\in \CC_D(h)$
such that $w_i\neq a_i$ for every~$i$.  If the $a_i$ are chosen
independently at random from~$\Pr^1(k)$, the expected number of
such words~$w$ is $\bigl( q/(q+1) \bigr)^N \! \cdot \, \# (\CC_D(h))$.
These words constitute a code of length~$N$\/
and minimal distance $\geqs\; N-2h$\/ over an alphabet of size~$q$.
But by (\ref{eq:Masymp}) the size of this code is
within a subexponential factor $\exp(o(N))$ of $q^{2h-g}$,
the Riemann-Roch lower bound on the number of words in the Goppa code
with the same alphabet size, length, and designed distance!
Since an average degradation of $\CC_D(h)$ is thus
asymptotically as good as a Goppa code,
we may justifiably claim that $\CC_D(h)$ itself is better than Goppa.

\section{Proofs}

We establish the lower bound $N-2h$\/
on the minimal distance of $\CC_D(h)$,
the independence of $\CC_D(h)$ of the choice of $\vph_i$,
and the isomorphism $\CC_D(h) \cong \CC_{D'}(h)$ when the
degree-zero divisors $D,D'$ are linearly equivalent.
We then prove the asymptotic formula~(\ref{eq:Masymp}) for $M(h,C)$,
and indicate how to modify our analysis to estimate the size
of individual codes $\CC_D(h)$.

\subsection{The distance bound}

\begin{prop}
\label{prop:distance}
Let $D$\/ be a divisor of degree~$0$ on a curve $C$\/ over~$k$,
and suppose $f_1,f_2$ are distinct sections of~$L_D$
of degrees $h_1,h_2$.  Then the words associated to $f_1,f_2$
by~(\ref{eq:NtupleD}) agree on at most $h_1+h_2$ coordinates.
In particular, $\CC_D(h)$ has minimal distance at least $N-2h$.
\end{prop}

\begin{proof}
We may assume that the $f_j$ are nonzero.
Let $E_1,E_2$ be the divisors $(f_1)+D$, $(f_2)+D$.
These $E_j$ are degree-$0$ divisors whose positive and negative parts
$E_j^{+},E_j^{-}$ each have degree $h_j$.  Set $f=f_1-f_2$,
a nonzero rational function on~$C$.  If $f_1,f_2$ agree on
the $i$-th coordinate then $P_i$ is either a pole of both
$\vph_i f\0_1$ and $\vph_i f\0_2$ or a zero of $\vph_i f$.  Let
\begin{equation}
\label{eq:Sdef}
S := \{ i : 1 \;\leqs\; i \;\leqs\; N, \;
(\vph_i f\0_1)(P_i) = (\vph_i f\0_2)(P_i) = \infty \},
\end{equation}
and $m=\#(S)$.
Then the negative part of the degree-zero divisor $D+(f)$
is bounded above by
$E_1^- + E_2^- - \sum_{i\in S} (P_i)$,
and thus has degree at most $h_1+h_2-m$.
Thus the positive part of $D+(f)$ also has degree at most $h_1+h_2-m$.
Hence there are at most $h_1+h_2-m$ choices of~$i$ for which
$(\vph_i f)(P_i)=0$.  Since there are $m$ common poles, we deduce
that the words associated to $f_1,f_2$ have at most
$(h_1+h_2-m)+m = h_1+h_2$ common coordinates, as claimed.
\end{proof}

\subsection{Easy isomorphisms}

\begin{lemma}
\label{lem:indep}
All choices of $\vph_i$ in~(\ref{eq:NtupleD}) yield equivalent codes.
\end{lemma}

\begin{proof}
Let $\psi_i$ be any other choice, and set $\theta_i = \psi_i/\vph_i$.
Then $\theta_i$ is a rational function on~$C$\/
with neither pole nor zero at~$P_i$.
Thus using $\psi_i$ instead of $\vph_i$ in~(\ref{eq:NtupleD})
multiplies the $i$-th coordinate of every word
by the nonzero scalar $\theta_i(P_i)$, for each~$i$.
Since each coordinate is changed
by a permutation of the alphabet $k\cup\{\infty\}$,
an equivalent code results.
\end{proof}

\begin{lemma}
\label{lem:lineq}
If $D,D'$ are linearly equivalent divisors of degree~$0$
then the codes $\CC_D(h)$, $\CC_{D'}(h)$ are isomorphic.
\end{lemma}

\begin{proof}
Let $D'-D$\/ be the divisor of the function~$g$.
Then $f$\/ is a rational section of degree $\leqs\; h$\/ of~$D'$
if and only if $fg$ is a rational section of degree $\leqs\;h$\/ of~$D$.
This identifies $\CC_D(h)$ and $\CC_{D'}(h)$ as sets.  Having chosen
$\vph_i$ for~$D$, we may choose $\varphi'_i := g \vph_i$ for~$D'$.
Then (\ref{eq:NtupleD}) gives the same coordinates for~$f$\/
as an element of $\CC_{D'}(h)$
that $fg$ has as an element of $\CC_D(h)$.
This identifies $\CC_D(h)$ and $\CC_{D'}(h)$ as error-correcting codes.
\end{proof}

Some {\sc remarks} on automorphisms: for nonzero $\theta\in k$\/
we have an isomorphism $f\mapsto \theta f$\/ from $\CC_D(h)$ to itself.
Thus the multiplicative group $k^*$ acts on~$\CC_D(h)$.
For general $C,D,h$\/ we expect that this is the full automorphism group
of $\CC_D(h)$.  By comparison, the Goppa codes, being linear, have
many more automorphisms: translation by any codeword, as well as
scalar multiplication.  Like the Goppa codes, our $\CC_D(h)$ can inherit
more symmetries from automorphisms of~$C$\/ and/or~$k$.
Thus if~$C$\/ has an automorphism
taking~$D$\/ to a divisor linearly equivalent to~$D$\/
then $\CC_D(h)$ inherits this automorphism by Lemma~\ref{lem:lineq}.
In particular, every automorphism of~$C$\/ acts in $\CC_0(h)$.
Likewise, if $C,D$\/ can be defined over a subfield~$k_0$ of~$k$\/
then Gal$(k/k_0)$ acts on $\CC_D(h)$.
Finally, $\CC_0(h)$ also has automorphisms by the group $\PGL_2(k)$,
which acts on~$\Pr^1(k)$ by fractional linear transformations.
Indeed, each $\gamma\in\PGL_2(k)$ yields the automorphism
$f\mapsto \gamma\circ f$\/ of~$\CC_0(h)$.
These automorphisms have no Goppa-code analogue.

\subsection{The average size {\large $M(h,C)$} of {\large $\CC_D(h)$}}
This requires more work.  For instance, the functions in $\CC_0(h)$
can be regarded the elements of height $\leqs\; h$\/
of the function field $k(C)$.
By a function-field analogue of a theorem of Schanuel~\cite{Schanuel},
announced by Serre~\cite[p.19]{Serre}
and proved by DiPippo~\cite{DiPippo} and Wan~\cite{Wan}
(independently but in the same way),
for any genus-$g$ curve~$C$\/ over~$k$\/ the number of such elements
is asymptotic to
\begin{equation}
\label{eq:Schanuel}
q^{2h+1-g} \, \frac{L_C(1)}{L_C(2)}
\end{equation}
as $h\ra\infty$, where $L_C$ is the $L$-function of the curve
(defined below).  We shall see later that
\begin{equation}
\label{eq:L1/L2}
\frac{L_C(1)}{L_C(2)} = \bigl( (q+1)/q \bigr)^{N+o(N)}
\end{equation}
if $C$\/ is an asymptotically optimal curve.
The same formula can be obtained for the number
of rational sections of~$L_D$ of degree at most~$h$.
But we need formulas valid not for $h\ra\infty$ but for $h<N/2$,
and this requires explicit and sufficiently small error terms
in the asymptotic formula~(\ref{eq:Schanuel}).

It is enough to count the elements of $\CC_D(h) - \CC_D(h-1)$,
which are rational sections~$f$\/ of~$L_D$ of degree exactly~$h$.
These are the functions whose divisors are of the form $E^+ - E^- - D$
where $E^+,E^-$ are effective divisors of degree exactly~$h$
with disjoint supports.
Necessarily $E^+ - E^-$ is linearly equivalent to~$D$.
Conversely, for each ordered pair $(E^+,E^-)$
of \hbox{degree-$h$} effective divisors with disjoint supports
such that $E^+ - E^- \sim D$,
there are $q-1$ rational functions~$f$\/
whose divisor is $E^+ - E^- - D$.
Thus $\#\bigl( \CC_D(h) - \CC_D(h-1) \bigr)$
is $(q-1)$ times the number of such ordered pairs $(E^+,E^-)$.
Averaging over~$D$\/ in~$J_C$ lets us ignore the condition
$E^+ - E^- \sim D$.

Now it is easy to count pairs $(D^+,D^-)$ of effective divisors
of degree~$n$ without the additional condition of disjoint supports:
the count is $M_n^2$,
where $M_n$ is the number of effective divisors of degree~$n$.
But each such pair $(D^+,D^-)$ is uniquely $(E+E^+,E+E^-)$
for some effective divisors $E,E^+,E^-$
with the supports of $E^+,E^-$ disjoint.
Thus
\begin{equation}
\label{eq:convol}
M_n^2 = \sum_{h=0}^n M_{n-h} A_h,
\end{equation}
where $A_0=1$, and for $h=1,2,3,\ldots$ we define
\begin{equation}
A_h := \frac1{q-1} \#(J_C) \bigl( M(h,C) - M(h-1,C) \bigr),
\label{eq:Ah}
\end{equation}
which is the number of pairs $(E^+,E^-)$
of effective divisors of degree~$h$ and disjoint supports.
The identity (\ref{eq:convol}) states that
the sequence $\{ M_n^2 \}$ is the convolution
of $\{M_n\}$ with $\{A_h\}$.  Thus
\begin{equation}
\label{eq:MM=M*A}
\sum_{h=0}^\infty A_h z^h =
\sum_{n=0}^\infty M_n^2 z^n \left/
\sum_{n=0}^\infty M_n z^n \right.
= Z_2(z) \, \big/ \, Z_1(z),
\end{equation}
where
\begin{equation}
\label{eq:Zm}
Z_m(z) := \sum_{n=0}^\infty M_n^m z^n.
\end{equation}
This leads us to study the functions $Z_1(z),Z_2(z)$.

Now $Z_1(z)$ is closely related to
the {\em zeta function} $\zeta_C$ of~$C$, defined by
\begin{equation}
\label{eq:zeta}
\zetaC(s) := \sum_{n=0}^\infty M_n q^{-ns}.
\end{equation}
Indeed $\zetaC(s) = Z_1(q^{-s})$.  Define
\begin{equation}
\label{eq:L}
L_C(s) := (1-q^{-s}) \, (1-q^{1-s}) \, \zetaC(s).
\end{equation}
It is known that $L_C(s)$, the {\em L-function} of~$C$,
is a polynomial of degree~$2g$ in~$q^{-s}$, of the form
\begin{equation}
\label{eq:lambda}
L_C(s)=\prod_{j=1}^{2g} (1- \lambda_j q^{-s}),
\end{equation}
where the $\lambda_j$, the ``eigenvalues of Frobenius'' for~$C$,
are $g$ conjugate pairs of complex numbers,
all of absolute value $q^{1/2}$.
(This is the ``Riemann hypothesis'' for~$L_C$,
here a celebrated theorem of Weil.)  Hence
\begin{equation}
\label{eq:Z1}
Z_1(z) = \prod_{j=1}^{2g} (1-\lambda_j z)
            \bigl/ \bigl( (1-z)(1-qz) \bigr) \, .
\end{equation}
This yields the exact formula
\begin{equation}
\label{eq:Mlarge}
M_n = \frac{q}{q-1} L_C(1) \, (q^n-q^{g-1})
\quad {\rm for}\ n > 2g-2.
\end{equation}
It follows that $Z_2(z)$ has a simple pole at $z=q^{-2}$
with residue
\begin{equation}
\label{eq:residue1}
- \left( \frac{q}{q-1} L_C(1) \right)^2
\end{equation}
and no other singularities
except for simple poles at $z=q^{-1}$ and $z=1$.
Thus $Z_2(z)/Z_1(z)$ has a simple pole at $z=q^{-2}$ with residue
\begin{equation}
\label{eq:residue}
- \frac
  {( \frac{q}{q-1} L_C(1) )^2}
  {(1-q^{-2})^{-1}\9 (1-q^{-1})^{-1}\9 L_C(2)}
= - \frac{q+1}{q} \, \frac{L_C^2(1)} {L_C(2)} \, ,
\end{equation}
and no other poles with $|z|<q^{-1/2}$, whence
\begin{equation}
\label{eq:Ahasymp}
A_h = \frac{q+1}{q} \, \frac{L_C^2(1)}{L_C(2)} q^{2h}
+ O_\epsilon (q^{(\frac12+\epsilon)h})
\end{equation}
as $h\ra\infty$.

It is further known that $\#(J_C)$ is given by the formula
\begin{equation}
\label{eq:|J|}
\#(J_C) = q^g L_C(1) = \prod_{j=1}^{2g} (1-\lambda_j)
\end{equation}
(``Dirichlet class number formula'' for function fields).
Hence
\begin{equation}
\label{eq:Mhasymp}
\frac{\sum_{n=0}^h A_n}{\#(J_C)} =
\frac{q}{q-1} \, \frac{L_C(1)}{L_C(2)} \, q^{2h-g}
 + O_\epsilon (q^{(\frac12+\epsilon)h}),
\end{equation}
so we have recovered (\ref{eq:Schanuel}) averaged over~$J_C$.
Still, we need estimates on $A_0+\ldots+A_h$ for $h<N/2$,
not as $h\ra\infty$.

To go further we use the distribution of the $\lambda_j$
on the circle $|\lambda|^2=q$.  Let $\alpha_j\in\R/2\pi\Z$
be the argument of $\lambda_j$:
\begin{equation}
\label{eq:alpha}
\lambda_j = q^{1/2} e^{i\alpha_j}.
\end{equation}
It is known that a family of curves $C$\/ is asymptotically optimal
if and only if
\begin{equation}
\label{eq:alpha_fourier}
\frac1{2g} \sum_{j=1}^{2g} e^{i r \alpha_j} \ra
\frac{1-\sqrt{q}}{2} \, q^{-|r|/2}
\end{equation}
for each nonzero integer~$r$
(see for instance ``Remark~1'' in~\cite{DV}).
Thus if $C$\/ is asymptotically optimal
then for any continuous function $\phi: \R/2\pi\Z \ra \C$ we have
\begin{equation}
\label{eq:phisum}
\frac1{2g} \sum_{j=1}^{2g} \phi(\alpha_j) \ra
a_0 + \frac{1-\sqrt{q}}{2} \sum_{r=1}^\infty q^{-r/2} (a_r + a_{-r}),
\end{equation}
where the $a_r$ are the Fourier coefficients of~$\phi$:
\begin{equation}
\label{eq:phi_fourier}
\phi(\alpha) \sim \sum_{r=-\infty}^\infty a_r e^{i r \alpha}.
\end{equation}
Since $\log L_C(s)  = \sum_{j=1}^{2g} \log(1-\lambda_j q^{-s})$ and
\begin{equation}
\label{eq:log_fourier}
\log(1 - z q^{1/2} e^{i\alpha}) =
- \sum_{r=1}^\infty \frac{(q^{1/2}z)^r}{r} e^{i r \alpha}
\end{equation}
for $|z|<q^{-1/2}$, we calculate
\begin{eqnarray}
\frac1g \log L_C(s) & \ra &
(q^{1/2}-1) \sum_{r=1}^\infty \frac{q^{-rs}}{r}
\nonumber \\
& = & -(q^{1/2}-1) \log(1-q^{-s})
\label{eq:log_L}
\end{eqnarray}
for all $s$ in $\{\Re(s)>1/2\}$, uniformly in any half-plane
$\Re(s) > \sigma$ with $\sigma>1/2$.  In particular,
since $N/g \ra q^{1/2}-1$ for our curves, we have
\begin{equation}
\label{eq:log_Lquot}
\frac1N \log \frac{L_C(1)}{L_C(2)} \ra
-\log \frac{1-q^{-1}}{1-q^{-2}} = \log \frac{q+1}{q},
\end{equation}
as we claimed in (\ref{eq:L1/L2}).

We can now prove:

\begin{theorem}
\label{thm:rho0}
For $\rho>0$ define $B(\rho)$ by
\begin{equation}
\label{eq:B}
B(\rho) := \min_{q^{-2} \leqs r \leqs q^{-3/2}}
 r^{-\rho/2} \frac{1-r}{(1-q^{-1})(1-qr)}.
\end{equation}
Then
\begin{equation}
\label{eq:Berror}
A_h = \frac{q+1}{q} \, \frac{L_C^2(1)}{L_C(2)} q^{2h}
+ O\Bigl( B\bigl(\frac{2h}{N}\bigr)^N\9 \exp \,o(N) \Bigr).
\end{equation}
We have
\begin{equation}
B(\rho) \ \leqs\ q^\rho (q+1)/(q-1)
\label{eq:small_B}
\end{equation}
for all $\rho>0$, with strict inequality if $\rho > 2q/(q^2-1)$.
If $C$\/ is asymptotically optimal (i.e., if $C$\/ varies
in a family of curves of genus \hbox{$g\ra\infty$}
with $N\sim(q^{1/2}-1)g$ rational points),
and for each~$C$\/ we choose $h$ with $\inf(h/N) > q/(q^2-1)$,
then $\log M(h,C)$ is given asymptotically by~(\ref{eq:Masymp}).
\end{theorem}

\begin{proof}

We estimate the error in (\ref{eq:Ahasymp}) using contour integration.
By (\ref{eq:MM=M*A}) and the discussion around (\ref{eq:residue})
we have
\begin{equation}
\label{eq:Aherror}
A_h - \frac{q+1}{q} \, \frac{L_C^2(1)}{L_C(2)} q^{2h} 
= \frac1{2\pi i}
\oint_{|z|=r} \frac{Z_2(z)}{Z_1(z)} \, \frac{dz}{z^{h+1}}
\end{equation}
for any $r\in(q^{-2},q^{-3/2})$.
(In fact we obtain (\ref{eq:Aherror}) for all $r\in(q^{-2},q^{-1/2})$,
but we shall soon need to assume $r<q^{-3/2}$.)
On the circle $|z|=r$ we have
\begin{equation}
\label{eq:|Z1|}
\log |Z_1(z)| = -N \log |1-z| + o(N)
\end{equation}
by~(\ref{eq:log_L}).
We estimate $|Z_2(z)|$ by using another contour integral
to express $Z_2(z)$ in terms of $Z_1$:
\begin{lemma}
\label{lem:Z2}
For all $z\neq q^{-1}$ with $q^{-2} < |z| < 1$ we have
\begin{eqnarray}
Z_2(z) & = & \frac1{2\pi i}
   \oint_{|w|=|z|^{1/2}} Z_1(w) Z_1\bigl(\frac{z}{w}\bigr) \frac{dw}{w}
\nonumber \\
      & + & 2 \, \frac{q^2}{q-1} L_C(1) Z_1(qz).
\label{eq:Z2int}
\end{eqnarray}
\end{lemma}

\begin{proof}
Consider first $z$ with $0<|z|<q^{-2}$.
For such~$z$ we obtain
\begin{equation}
\label{eq:Z2small}
Z_2(z)  =  \frac1{2\pi i} \oint_{|w|=|z|^{1/2}}
            Z_1(w) Z_1\bigl(\frac{z}{w}\bigr) \frac{dw}{w}
\end{equation}
by integrating termwise the product
of the absolutely convergent series~(\ref{eq:Zm})
for $Z_1(w)$ and $Z_1(z/w)$.  For any $z$ other than
$0,q^{-2},q^{-1},1$, the integrand extends to a meromorphic
function on~$\C$ with simple poles at $w=z,qz,1/q,1$
and a multiple pole at $w=0$.
The contour in~(\ref{eq:Z2small}) encloses the poles $0,z,qz$
but not the poles $1/q,1$.
Thus analytic continuation gives
\begin{equation}
\label{eq:Z2cont}
Z_2(z) = \frac1{2\pi i} \oint Z_1(w)
          Z_1\bigl(\frac{z}{w}\bigr) \frac{dw}{w}
\end{equation}
for all $z\notin\{q^{-2},1\}$, for any contour that encloses $0,z,qz$
but not $1/q,1$.
Now when $q^{-2}<|z|<1$ the contour in~(\ref{eq:Z2int})
encloses $0,z,1/q$ but not $qz,1$.
Thus we can evaluate the contour integral in~(\ref{eq:Z2int})
by starting from (\ref{eq:Z2cont}),
adding the residue at $1/q$, and subtracting the residue at $qz$.
The former residue is $-(q^2/(q-1)) \, L_C(1) Z_1(qz)$,
and the latter is $+(q^2/(q-1)) \, L_C(1) Z_1(qz)$.
This proves (\ref{eq:Z2int}).
\end{proof}

Thus (\ref{eq:Aherror}) is
\begin{equation*}
\frac1{2\pi i}
\oint_{|z|=r} \frac1{Z_1(z)}
 \left(\frac{2q^2}{q-1} L_C(1) Z_1(qz) \right.
\qquad\qquad
\end{equation*}
\vspace*{-1ex}
\begin{equation}
\qquad\qquad
 + \frac1{2\pi i} \left.
 \oint_{|w|=|z|^{1/2}} Z_1(w) Z_1\bigl(\frac{z}{w}\bigr) \frac{dw}{w}
 \right) \, \frac{dz}{z^{h+1}}.
\label{eq:Ahiint}
\end{equation}
We use (\ref{eq:log_L},\ref{eq:|Z1|}) to estimate both parts of this.
For the single integral, we find
\begin{eqnarray}
&&\!\!\!\!
\log
\left|\frac{2q^2}{q-1} \frac{Z_1(qz)}{Z_1(z)} \frac{L_C(1)}{z^h} \right|
\nonumber\\
\label{eq:int<}
&&\!\!\!\!
 = -h \log r + N \log \frac{|1-z|}{(1-q^{-1})|1-qz|} + o(N)
\nonumber\\
&&\!\!\!\!
\leqs -h \log r + N \log \frac{1-r}{(1-q^{-1})(1-qr)} + o(N).\quad
\end{eqnarray}
Thus the single integral is $O(B(2h/N)^N \exp o(N))$.
We shall show that the double integral is exponentially smaller
than $B(2h/N)^N$; this will prove (\ref{eq:Berror}).
To estimate the integrand, let $w'=z/w$, so $z=ww'$ and
\begin{eqnarray}
&&\!\!\!\!
\log \left|
 \frac{2q^2}{q-1} \, \frac{Z_1(w)Z_1(w')}{Z_1(ww')} \, \frac1{z^h}
\right|
\nonumber\\
&&\!\!\!\!
= -h \log r + N \log \left|\frac{(1-w)(1-w')}{1-ww'}\right| + o(N).
\label{eq:iint<}
\end{eqnarray}
Here $|w|=|w'|=r^{1/2}$, so
\begin{eqnarray}
\left|\frac{(1-w)(1-w')}{1-ww'}\right|
&=& \left| 1 + \frac{w}{(1-w)} + \frac{w'}{(1-w')} \right|
\nonumber\\
\leqs \ 1 + 2\frac{r^{1/2}}{1-r^{1/2}}
&=& \frac{1+r^{1/2}} {1-r^{1/2}}.
\label{eq:w,w'}
\end{eqnarray}
Thus our proof of (\ref{eq:Berror}) will be complete once we show
\begin{equation}
\label{eq:i>ii}
\frac{1-r}{(1-q^{-1})(1-qr)} \ \geqs\ \frac{1+r^{1/2}} {1-r^{1/2}},
\end{equation}
or equivalently
\begin{equation}
\label{eq:i>ii.1}
(1-r^{1/2})^2 \;\geqs\; (1-q^{-1}) (1-qr);
\end{equation}
and this follows from the observation that
\begin{equation}
\label{eq:i>ii.2}
(1-r^{1/2})^2 - (1-q^{-1}) (1-qr) = q (r^{1/2}-q^{-1})^2.
\end{equation}

It remains to prove (\ref{eq:small_B}) and to show that
the ``main term'' in (\ref{eq:Berror}) is indeed exponentially larger
than the ``error term'' as long as $\inf(2h/N) > 2q/(q^2-1)$.
By (\ref{eq:log_L}), the main term is
\begin{equation}
\label{eq:mainterm}
q^{2h} \left(\frac{q+1}{q-1}\right)^N \exp(o(N)).
\end{equation}
Thus strict inequality in the upper bound (\ref{eq:small_B})
is what we need to show that
(\ref{eq:mainterm}) exceeds the ``error term''$\!$.
The ratio between $B(\rho)$ and the claimed upper bound is
\begin{equation}
\label{eq:Bratio}
q^{-\rho} \, \frac{q-1}{q+1} \, B(\rho) =
\frac{q}{q+1} \min_{q^{-2} \leq r \leq q^{-3/2}}
 (q^2 r)^{-\rho/2} \frac{1-r}{1-qr}.
\end{equation}
Trying $r=q^{-2}$ we find that
\begin{equation}
\label{eq:Bratio<=1}
q^{-\rho} \, \frac{q-1}{q+1} \, B(\rho)
\ \leqs\
\frac{q}{q+1} \, \frac{1-q^{-2}}{1-q^{-1}} = 1,
\end{equation}
so the upper bound holds for all~$\rho$.  Moreover the bound is strict
if $r^{-\rho/2} (1-r)/(1-qr)$ is a decreasing function of~$r$
at $r=q^{-2}$.  We calculate that the logarithmic derivative of
$r^{-\rho/2} (1-r)/(1-qr)$ at $r=q^{-2}$ is
\begin{equation}
\label{eq:dlog}
-\frac{q^2}{q^2-1} \bigl( (q^2-1) \frac\rho2 - q \bigr).
\end{equation}
This is negative once $\rho > 2q/(q^2-1)$, so Theorem~\ref{thm:rho0}
is proved.~~\lower1pt\hbox{\large{$\Box\Box$}}
\end{proof}

\subsection{The size of individual codes {\large $\CC_D(h)$}}

We showed above that $\#(\CC_D(h)-\CC_D(h-1))$
is $(q-1)$ times the number of ordered pairs $(E^+,E^-)$
of effective \hbox{degree-$h$}\/ divisors with disjoint supports
such that \hbox{$E^+-E^- \sim D$.}  Call this number $A_h(D)$,
so that the total count $A_h$ introduced in~(\ref{eq:Ah})
is $\sum_{D\in J_C} A_h(D)$.  We expect that $A_h(D)$
is approximated by $A_h/\#(J_C)$ if $h$\/ is large enough.

To prove this we use a known device from analytic number theory:
for each character~$\chi$ of the finite abelian group $J_C$, define
\begin{equation}
\label{eq:Ahchi}
A_h(\chi) := \sum_{D\in J_C} \chi(D) A_h(D).
\end{equation}
This is the sum of $\chi(E^+ - E^-)$ over all ordered pairs
of effective divisors $E^+,E^-$ of degree~$h$ with disjoint supports.
{}From the $A_h(\chi)$ we can recover $A_h(D)$ by the usual formula
\begin{equation}
\label{eq:chi_to_D}
A_h(D) = \frac1{\#(J_C)} \sum_\chi \overline\chi(D) A_h(\chi).
\end{equation}
When $\chi$ is the trivial character
(the character sending all of~$J_C$ to~$1$),
the sum $A_h(\chi)$ reduces to $A_h$;
we expect that the other $A_h(\chi)$ will be smaller.
As with $A_h$, we analyze the $A_h(\chi)$ by comparing them with
\begin{equation}
\label{eq:Nnchi}
N_n(\chi) := \sum_{\deg(D^+)=\deg(D^-)=n} \chi(D^+ - D^-),
\end{equation}
the sum extending over all pairs of effective divisors $D^+,D^-$,
whether disjointly supported or not.  Again, any such pair is uniquely
$(E+E^+,E+E^-)$ with $E,E^+,E^-$ effective divisors such that
$E^+,E^-$ have disjoint supports; and necessarily
$\chi(E^+ - E^-) = \chi(D^+ - D^-)$.
Thus we have a convolution formula
\begin{equation}
\label{eq:convolchi}
N_n(\chi) = \sum_{h=0}^n M_{n-h} A_h(\chi),
\end{equation}
generalizing (\ref{eq:convol}).  We deduce that
\begin{equation}
\label{eq:Ahchiquot}
\sum_{h=0}^\infty A_h(\chi) z^h = Z_2(z,\chi)/Z_1(z),
\end{equation}
with $Z_1(z)=\sum_{n=0}^\infty M_n z^n$ as above and
\begin{equation}
\label{eq:Z2chi}
Z_2(z,\chi) := \sum_{n=0}^\infty N_n(\chi) z^n.
\end{equation}

We can factor $N_n(\chi)$ by writing
\begin{equation}
\chi(D^+-D^-) = \chi(D^+) \, \overline\chi(D^-).
\label{eq:chichibar}
\end{equation}
Since $D^\pm$ are not in general divisors of degree zero,
this requires that $\chi$ be extended from $J_C$ to the group
$\Pic(C)$ of linear equivalence classes of divisors on~$C$\/
of arbitrary degree.  For each $\chi$, choose an arbitrary extension
of~$\chi$ to a homomorphism from $\Pic(C)$ to the unit circle.
[For instance, fix a divisor $D_1$ of degree~$1$,
and let $\chi(D_1)$ be an arbitrary complex number of norm~$1$;
any such choice of $\chi(D_1)$ yields a unique extension
of~$\chi$ to $\Pic(C)$.]  Then
\begin{equation}
\label{eq:Nnchifactor}
N_n(\chi) = M_n(\chi) M_n(\overline\chi),
\end{equation}
where $M_n(\chi)$ is the sum of the values of~$\chi$
on effective divisors of degree~$n$.
[Changing $\chi(D_1)$ to $\beta\chi(D_1)$, for some $\beta\in\C$
of norm~$1$, multiplies $M_n(\chi)$ and $M_n(\overline\chi)$
by $\beta^n$ and $\beta^{-n}$ respectively,
and thus does not change their product.]

For a nontrivial character~$\chi$ we have $M_n(\chi)=0$
for all \hbox{$n>2g-2$,} because by Riemann-Roch
each degree-$n$ class in $\Pic(C)$ is represented
the same number of times in the sum $M_n(\chi)$.\footnote{
  This already suffices to show that as $h\ra\infty$ the formula
  \begin{equation}
  \label{AhDasymp}
  A_h(D) = \frac{q+1}{q} \, \frac{L_C(1)}{L_C(2)} q^{2h-g}
  + O_\epsilon (q^{(\frac12+\epsilon)h})
  \end{equation}
  holds not only on average over~$D$\/ (this average estimate
  is~(\ref{eq:Ahasymp})) but also for each~$D$.  We thus recover
  Schanuel's theorem with a sharp error term.  But again
  our present application requires estimates for $h\ll N$,
  not $h\ra\infty$.
  }
Thus
\begin{equation}
\label{eq:Lchi}
L(s,\chi) := \sum_n M_n(\chi) q^{-ns}
\end{equation}
is a finite sum.
This sum, called the \hbox{$L$-function} associated to~$\chi$,
is again known to satisfy a Riemann hypothesis, which yields
a factorization
\begin{equation}
\label{eq:Lchifactor}
\sum_{n=0}^{2g-2} M_n(\chi) z^n =
\prod_{j=1}^{2g-2} (1-\lambda_j(\chi)z)
\end{equation}
for some $\lambda_j(\chi)$ all of absolute value~$q^{1/2}$.
Unlike the eigenvalues of Frobenius~$\lambda_j$ for~$C$,
the $\lambda_j(\chi)$ are of unknown distribution
even for an asymptotically optimal~$C$.
Thus instead of asymptotic formulas for
\begin{equation}
\label{eq:Z1chi}
Z_1(z,\chi) := \sum_{n=0}^{2g-2} M_n(\chi) z^n
\end{equation}
we get only an upper bound:
\begin{equation}
\label{eq:Z1chibd}
|Z_1(z,\chi)| \ \leqs\ (1+q^{1/2}|z|)^{2g-2}
\end{equation}
for all $z\in\C$.  But an upper bound is all we need because
$|Z_1(z,\chi)|$ contributes only to the error terms $A_h(\chi)$,
$A_h(\overline\chi)$.  Since $Z_1(z,\chi)$ is a polynomial,
we need not worry about nonzero poles in the contour integral
\begin{equation}
\label{eq:Z2oint}
Z_2(z,\chi)  =  \frac1{2\pi i} \oint_{|w|=|z|^{1/2}}
            Z_1(w,\chi) Z_1\bigl(\frac{z}{w},\overline\chi\bigr)
            \frac{dw}{w}
\end{equation}
for $Z_2(z,\chi)$, which holds for all $z\neq 0$.  Therefore
\begin{equation}
\label{eq:Z2bd}
|Z_2(z,\chi)| < (1+\sqrt{q|z|}\,)^{4g-4}.
\end{equation}
Using contour integration about a circle of radius~$r$
to isolate the $z^h$ term of (\ref{eq:Ahchiquot}), we obtain
\begin{equation}
\label{eq:Ahchibd}
|A_h(\chi)| < r^{-h} (1+\sqrt{qr}\,)^{4g} (1+r)^{N+o(N)}
\end{equation}
for any positive $r<q^{-1/2}$.  Minimizing this over~$r$,
summing over the $\#(J_C)$ choices of~$\chi$, and using
our known estimates for $A_h$ and $L_C(1)$, we find:
\begin{theorem}
\label{thm:rho1}
For $\rho>0$ define $B_1(\rho)$ by
\begin{equation}
\label{eq:B1}
B_1(\rho) := \frac{q+1}{q} q^\kappa
  \min_{r \leqs q^{-1/2}}
  r^{-\rho/2} (1+r) (1+\sqrt{qr}\,)^{4\kappa} ,
\end{equation}
where $\kappa := 1/(\sqrt{q}-1)$.  Then
\begin{equation}
\label{eq:B1error}
A_h(D) = \frac1{\#(J_C)} \left( A_h
  + O\Bigl( B_1\bigl(\frac{2h}{N}\bigr)^{N}\9 \exp \,o(N) \Bigr)
  \right)
\end{equation}
for every degree-$0$ divisor~$D$.
There exists a unique $\rho\0_1=\rho\0_1(q)>0$ such that
\begin{equation}
\label{eq:rho_1}
B_1(\rho\0_1) = q^{\rho\0_1} \frac{q+1}{q-1};
\end{equation}
$B_1(\rho) < q^\rho (q+1)/(q-1)$ for all $\rho>\rho\0_1$.
If $C$\/ is asymptotically optimal,
and for each~$C$\/ we choose $h$ with $\inf(2h/N) > \rho\0_1$,
then $\log \#(\CC_D(h))$ is given asymptotically by
\begin{equation}
\label{eq:Casymp}
\#(\CC_D(h)) = \left(\frac{q+1}{q}\right)^{N+o(N)} q^{2h-g}.
\end{equation}
\end{theorem}

\begin{proof}
Estimate (\ref{eq:B1error}) follows from (\ref{eq:chi_to_D})
and the bound (\ref{eq:Ahchibd}) on each term with $\chi$ nontrivial,
together with the facts $g/N \ra \kappa$ and
\begin{equation}
\label{eq:log_J}
\#(J_C) = q^g \left( \frac{q+1}{q} \right)^{N+o(N)}
\end{equation}
(see (\ref{eq:|J|},\ref{eq:log_L})).  For the remainder term
to be exponentially smaller we must have $h/N > q/(q^2-1)$
(from Thm.~\ref{thm:rho0}) and
\begin{equation}
\label{eq:small1}
B_1\bigl(\frac{2h}{N}\bigr) < q^{2h/N} \frac{q+1}{q-1}.
\end{equation}
The ratio between the two sides is
\begin{equation}
\label{eq:B1ratio}
q^\kappa \frac{q-1}{q} \min_{r \leqs q^{-1/2}}
  (q^2 r)^{-\rho/2} (1+r) (1+\sqrt{qr}\,)^{4\kappa} ,
\end{equation}
where again $\rho=2h/N$.  For all $r \;\leqs\; q^{-2}$,
the product (\ref{eq:B1ratio}) exceeds $(1-q^{-1})q^\kappa>1$.
For $r=q^{-1/2}$ the product clearly falls below~$1$
once $\rho$ is large enough.  Thus (\ref{eq:B1ratio}) equals~$1$
for some $\rho\0_1$, with the minimum attained at some $r>q^{-2}$;
since $(q^2 r)^{-\rho/2} (1+r) (1+\sqrt{qr}\,)^{4\kappa}$
is a decreasing function of~$\rho$ for that~$r$, the inequality
(\ref{eq:small1}) holds for all $2h/N>\rho\0_1$.
It is not hard to check that $\rho_1 > 2q/(q^2-1)$ ---
even the lower bound $q^\kappa \frac{q-1}{q} q^{-\frac32q/(q^2-1)}$
on~(\ref{eq:Bratio}) suffices for this.
The claim (\ref{eq:Casymp}) now follows
{}from (\ref{eq:B1}) and Thm.~\ref{thm:rho0}.
\end{proof}

The following short table lists $\rho\0_1$ rounded to four decimals
for $q=q_0^2$ and $q\0_0$ a prime power $\leqs\; 16$:

\noindent
\centerline{
\begin{tabular}{c|ccccc}
   $q$     &  $2^2$ &  $3^2$ &  $4^2$ &  $5^2$ &  $7^2$ \\ \hline
$\rho\0_1$ & 4.3461 & 1.8541 & 1.1606 & 0.8348 & 0.5276
\end{tabular}
}

\vspace*{2ex}

\centerline{
\noindent
\begin{tabular}{c|ccccc}
   $q$     &  $8^2$ &  $9^2$ & $11^2$ & $13^2$ & $16^2$ \\ \hline
$\rho\0_1$ & 0.4440 & 0.3827 & 0.2990 & 0.2448 & 0.1919
\end{tabular}
}

Since the definition of $\CC_D(h)$ requires $2h<N$,
we must have $\rho<1$, so the threshold $\rho\0_1$ is too high
for $q=4,9,16$.  For these small~$q$, we get information
only about the average size $M_h(C)$ of the codes $\CC_D(h)$
with small~$\delta$.
But it is only for $q \;\geqs\; 49$ that any of the
algebraic-geometry codes improve on \GV.  For $q=49$ it turns out
that $\rho\0_1$ is larger than the maximal $2h/N$\/ for which
$M_h(C)$ attains or exceeds the \GV\ bound.  For $q \;\geqs\; 64$,
we find that $\rho\0_1$ is within the range of codes whose
average size $M_h(C)$ improves on \GV; thus in each case we have
a subrange in which each individual code $\CC_D(h)$
is known to be exponentially larger than the \GV\ bound.
As $q$ increases, $\rho\0_1(q)\ra 0$, so this subrange of $2h/N$\/
values covers almost all of $(0,1)$.

\section{Problems}

\subsection{New problems in computational algebraic geometry}
A new construction of error-correcting codes
automatically raises new decoding problems.
When the codes come from algebraic curves,
these problems can be stated in terms of the geometry of the curves.
For example, for $\CC_0(h)$, the problem of nearest-neighbor decoding
is a special case of the following problem:

\vspace{3ex}

{\sc Problem 1.}
{\em
Given: an algebraic curve~$C$\/ of genus~$g$ over a field~$k$\/;
a list $(P_1,\ldots,P_N)$ of \hbox{$k$-rational} points of~$C$\/;
an \hbox{$N$\/-tuple} $(w_1,\ldots,w_N)$ in $(\Pr^1(k))^N$;
and integers $h,e\;\geqs\;0$.  Find a rational function~$f$\/
of degree at most~$h$\/ on~$C$\/ such that $f(P_i)=w_i$
for each~$i$ with at most $e$ exceptions,
assuming that at least one such~$f$\/ exists.
}

\vspace{3ex}

Similarly for $\CC_D(h)$:

\vspace{3ex}

{\sc Problem 1'.}
{\em
Given: an algebraic curve~$C$\/ of genus~$g$ over a field~$k$\/;
a divisor~$D$\/ of degree zero on~$C$\/; 
a list $(P_1,\ldots,P_N)$ of \hbox{$k$-rational} points of~$C$,
and functions $\vph_i$ whose divisor
has the same order at~$P_i$ as~$D$\/;
an $N$\/-tuple $(w_1,\ldots,w_N)$ in $(\Pr^1(k))^N$;
and integers $h,e\;\geqs\;0$.
Find a rational section~$f$\/ of~$D$\/
of degree at most~$h$\/ on~$C$\/ such that $(\vph_i f)(P_i)=w_i$
for each~$i$ with at most $e$ exceptions,
assuming that at least one such~$f$\/ exists.
}

\vspace{3ex}

By Prop.~\ref{prop:distance},
if $2(h+e)<N$\/ then $f$\/ is uniquely determined;
if $2(h+e)$ equals or exceeds $N$, but not by too much,
one might still hope that there are few enough spurious~$f$\/
that ``list decoding'' (that is, finding all possible~$f$, not just one)
may be feasible as in~\cite{GuSu,ShWa}.

The special case $e=0$ of Problem~1 or~1'
is the {\em error detection} or {\em recognition}\/ problem:
{\em is a given word in the code}?  For a Goppa code,
the recognition problem is readily solved
in time polynomial in the length of the code:
the code is linear, so recognition reduces to linear algebra.
But the new codes $\CC_D(h)$ are nonlinear,
and an efficient error-detection algorithm is not obvious.

Another, possibly even more fundamental, difficulty is
enumerating $\CC_D(h)$.  To use $\CC_D(h)$ in any error-correcting
application other than the highly unlikely application
of transmitting the values of a low-degree rational section of~$D$,
one must have an efficient means of generating the $m$-th codeword
as a function of~$m$, and of inverting this function to recover
the integer~$m$ transmitted.  For a linear code with a known basis,
enumeration is no harder than recognition, but again the problem
seems nontrivial for our nonlinear codes $\CC_D(h)$.  It is not
necessary to enumerate {\em every} codeword: if $M<\#(\CC_D(h))$,
an efficiently computable and invertible injection
{}from $[M] := \{1,\ldots,M\}$ to $\CC_D(h)$
would still let us use an \hbox{$M$-word} subcode
of $\CC_D(h)$ for error-resistant communication.
But $M$\/ must not be so much smaller than $\#(\CC_D(h))$
as to reduce the asymptotic transmission rate.  Thus we ask:

\vspace{3ex}

{\sc Problem 2.}
{\em
Find $M=\#(\CC_D(h))^{1-o(1)}$ and an injection
$\iota: [M] \hookrightarrow \CC_D(h)$
such that both $\iota$ and the inverse function
$\iota^{-1} : \iota([M]) \ra [M]$
are efficiently computable.
}

\vspace{2ex}

\subsection{Solutions for {\large $C$}\/ of genus zero}

We show that both Problems~1 and~2 have polynomial-time solutions
when $C$\/ has genus zero.  (In that case, all degree-zero divisors
are linearly equivalent, so Problems 1 and~1' are equivalent.)
This does not directly address the issue of using $\CC_D(h)$ for
error-resistant communications, because that application requires
curves of large genus; the most direct generalization of our solution
to arbitrary~$C$\/ requires exhaustion over~$J_C$ and thus takes time
exponential in the genus.  Nevertheless we have hope
that our solutions can be adapted to the large-genus case,
especially for Problems~1 and~1'.
This is because we solve Problem~1 in genus~zero
by adapting a known algorithm for decoding Reed-Solomon codes.
Goppa codes are large-genus generalizations of Reed-Solomon codes,
and can be decoded efficiently~\cite{GuSu,ShWa}.
It may be possible to combine ideas from these decoding algorithms
and our genus-zero solution of Problem~1
to solve that Problem in general.

In the genus-zero case,
all $\CC_D(h)$ with the same $q,h$\/ are isomorphic.
Thus we may and shall assume $D=0$, and call the codes simply
``$\CC(h)$'', suppressing the subscript.
This $\CC(h)$ consists of rational functions in one variable~$x$,
evaluated at $x=P_i$ (one of which may be $\infty$).
A rational function $f(\cdot)$ of degree~$h$\/
is a quotient $a(x)/b(x)$ of relatively prime polynomials
$a,b$ in~$x$ of degree $\leqs\; h$\/:
\begin{equation}
\label{eq:a,b}
a(x) = \sum_{j=0}^h a_j x^j,
\quad
b(x) = \sum_{j=0}^h b_j x^j,
\end{equation}
with the leading coefficients $a_h,b_h$ not both zero.
A condition $f(P_i)=w_i$ is a homogeneous linear equation
in the \hbox{$2h+2$} coefficients $a_j,b_j$.
(If $w_i=\infty$ the equation becomes \hbox{$b(P_i)=0$;}
if $P_i=\infty$ the equation is $a_h=w_i b_h$ if $w_i$ is finite,
$b_h=0$ if $w_i=\infty$.\footnote{
  As usual the special cases $P_i=\infty$, $w_i=\infty$
  that appear here and later can be avoided by using
  homogeneous coordinates on~$\Pr^1$ and regarding $f$\/
  as the quotient of two degree-$h$ homogeneous polynomials
  in two variables.
  })
Thus the recognition problem amounts to
solving the $N$\/ simultaneous linear equations coming from
$f(P_i)=w_i$, which we can do in time polynomial in~$N$.
We claim that every nonzero solution is proportional to $(a_j,b_j)$
and thus recovers the function $f=a/b$, as long as $2h<N$\/ ---
exactly the condition we imposed on~$h$\/ when we defined of $\CC(h)$.
Indeed, suppose $(a'_j,b'_j)$ is another solution,
yielding another rational function $f'=a'/b'$.
Then the polynomial $\Delta := a'b-ab'$, of degree at most~$2h$,
vanishes at all finite $P_i$,
and its $x^{2h}$ coefficient vanishes if some $P_i=\infty$.
Thus $\Delta$ is identically zero, and $f=f'$ as claimed.
If $f$\/ is of degree $<h$, the same argument shows that
the linear equations on $a_j,b_j$ will have a solution space
of dimension $h-\deg(f)+1$, and any nonzero solution vector
recovers~$f$\/ as~$a/b$.  We have thus solved the genus-zero case
of Problem~1 for $e=0$ and $2h<N$.

The same system of simultaneous linear equations
with $h$ replaced by $h+e$ also solves the genus-zero case of Problem~1
for any $e$ such that $2(h+e)<N$ --- that is,
for all $e$ less than half the designed distance $N-2h$\/ of the code.
To see this, suppose $f=a/b$ differs from the word~$w$ in at most
$e$ coordinates, and let $c(x)$ be an ``error-locating polynomial'':
a polynomial of degree at most~$e$ that vanishes at each finite $P_i$
where $f(P_i)\neq w_i$.  (If one of the errors is at $P_i=\infty$
then $c(x)$ has degree at most $e-1$.)  Then the coefficients of
the polynomials $ac$ and $bc$ satisfy the linear equations on the
coefficients of polynomials of degree $h+e$ whose quotient agrees
with~$w$ at all $P_i$.  Any solution $(a'_j,b'_j)$ of these equations
yields polynomials $a',b'$ such that $\Delta := c(a'b-ab')$, which now
is a polynomial of degree $\leqs\; 2(h+e)$, vanishes at all finite $P_i$
and has vanishing $x^{2(h+e)}$ coefficient if some $P_i=\infty$.
Again it follows that $\Delta=0$ identically and $f=a'/b'$.
Thus as claimed we can decode the codes $\CC(h)$
associated to $C=\Pr^1$ up to the error-correcting bound
$\frac12(N-1)-h$.

In the genus-zero case the enumeration problem also has
a polynomial-time solution, even without relaxing it
to a large subset of~$\CC(h)$ as in Problem~2.
When $C=\Pr^1$, the \hbox{$L$-function} of $C$\/ is the constant~$1$,
so we know $Z_1$ exactly, and thus also $Z_2$ and $A_h$.  We calculate:
\begin{equation}
\label{eq:P1(AB)}
Z_1(z) = \frac1{(1-z)(1-qz)}, \quad
M_n = \frac{q^{n+1}-1}{q-1},
\end{equation}
\vspace*{-1ex}
\begin{equation}
\label{eq:P1(CD)}
Z_2(z) = \frac{1+qz}{(1-z)(1-qz)(1-q^2 z)}, \quad
\frac{Z_2(z)}{Z_1(z)} = \frac{1+qz}{1-q^2 z},
\end{equation}
whence $A_h=q^{2h}+q^{2h-1}$ for $h>0$.  Since $A_0=1$,
\begin{equation}
\label{eq:P1(Sch)}
\#(\CC(h)) = 1 + (q-1)\sum_{i=0}^h A_i = q^{2h+1}
\end{equation}
(so the asymptotic formula (\ref{eq:Schanuel}) is exact here!\footnote{
  This result, but not the simpler proof we give next, already occurs
  in~\cite{DiPippo}, as a special case of a formula for~$\#(\CC_0(h)$
  depending only on the zeta function of~$C$\/ in the case that
  $C$ is hyperelliptic.
  })
We next construct a bijection~$\iota$ from $\CC(h)$
to a finite field $k'$ containing $k$ with degree $2h+1$.
Since $k'$ is readily enumerated
(choose a basis for $k'$ as a vector space over its prime field),
our bijection will yield a complete enumeration of $\CC(h)$.
To construct $\iota$, fix $x_0\in k'$ that generates~$k'$ over~$k$,
and define $\iota(f)=f(x_0)$ for all $f\in\CC(h)$.
Note that $f(x_0)$ cannot be $\infty$, because the denominator of~$f$\/
has degree at most $h<[k':k]$, and thus cannot vanish at~$x_0$.
Moreover, $\iota$ is an injection:
if $f_1,f_2$ are distinct rational functions of degree at most~$h$\/
we cannot have $f_1(x_0)=f_2(x_0)$, because then
$x_0$ would be a root of a polynomial of degree at most~$2h$,
and thus could not generate the field extension $k'/k$.
Since $\#(k')=\#(\CC(h))$ it follows that $\iota$ is a bijection.
To invert $\iota$, we must express any $x_1\in k'$ as $a(x_0)/b(x_0)$
for some polynomials $a,b$ of degrees $\leqs\; h$.
This, too, can be done by solving $2h+1$ simultaneous linear equations,
and thus in time polynomial in~$q$.  For instance, find the intersection
of the two $k$-vector subspaces
\begin{equation}
\label{eq:subspace1}
\{ a(x_0): a\in k[X], \ \deg(a) \;\leqs\; h \}
\end{equation}
and
\begin{equation}
\label{eq:subspace2}
\{ x_1 b (x_0): b\in k[X], \ \deg(b) \;\leqs\; h \},
\end{equation}
of dimension~$h+1$ in~$k'$.  Note that the intersection has dimension
at least $2(h+1)-(2h+1)=1$, and thus contains a nonzero vector.
This proves directly that the injection $\iota$ is onto, and thus
also completes an alternative proof of the formula (\ref{eq:P1(Sch)}).

{\sc Remark}\/: The algorithms in these section are polynomial-time
but far from optimal.  The simultaneous linear equations that arise
are of a special form that can be solved much more quickly by other
methods such as fast gcd's in~$k[X]$.

\subsection{Theoretical problems}

Our results also suggest at least three theoretical problems.
When $q \;\geqs\; 7^2$, it is known that Goppa's code can be modified
to improve on both \GV\ and (\ref{eq:Goppamax})
near the crossover points between these two lower bounds.

\vspace{3ex}

{\sc Problem 3.}
{\em
Does our construction of $\CC_D(h)$ admit similar improvements
near the crossover points between (\ref{eq:Newmax}) and the
\GV\ bound for codes over an alphabet of $q+1$ letters?
}

\vspace{3ex}

A second problem is whether the thresholds
$2q/(q^2-1)$ and $\rho\0_1(q)$
of Thms.~\ref{thm:rho0} and Thm.~\ref{thm:rho1} are best possible:

\vspace{3ex}

{\sc Problem 4.}
{\em
Can the bounds $2q/(q^2-1)$ and $\rho\0_1(q)$ be reduced?
In particular, can any of $\rho\0_1(4)$, $\rho\0_1(9)$, $\rho\0_1(16)$
be replaced by a threshold $<1$?
}

\vspace{3ex}

If $\rho\0_1(4)$ can be pushed below~$1$ then (\ref{eq:Masymp})
will yield a deterministic construction of
arbitrarily long algebraic-geometry codes over a five-letter alphabet
with $R,\delta$ both bounded away from zero.
Note that by~(\ref{eq:Goppamax}) Goppa codes do not do this when $q=4$.
For a five-letter alphabet,
Thm.~\ref{thm:rho0} proves the existence of such codes,
but does not let us specify one in time polynomial in~$N$,
because of the averaging over~$J_C$.  We may thus ask:

\vspace{3ex}

{\sc Problem 5.}
{\em
Is it possible to compute, in polynomial or random polynomial time,
a choice of~$D$\/ that makes $\CC_D(h)$ at least as large as average,
and thus with $R,\delta$ both provably bounded away from zero?
}

Finally, a more speculative kind of problem
concerns our earlier observation
that degrading $\CC_D(h)$ to a \hbox{$q$-letter} alphabet
yields nonlinear codes with exactly the same $R,\delta$ as Goppa codes.
Is this more than a coincidence?  That is,

\vspace{3ex}

{\sc Problem 6.}
{\em
Give a conceptual explanation for the factor $((q+1)/q)^N$
in~(\ref{eq:Masymp}), and for the fact that it exactly cancels
the degradation factor $(q/(q+1))^N$.
}

\section{Acknowledgments}
Thanks to the Packard foundation for financial support,
to Joel Rosenberg for a careful reading of an earlier draft,
and to Stephen DiPippo for the references~\cite{DiPippo,Wan}.

\balancecolumns


\begin{thebibliography}{00}

\bibitem{DV} V.~G. Drinfeld and S.~G. Vl\u{a}du\c{t}.
\newblock  The number of points of an algebraic curve.
\newblock  {\em Functional Anal.\ Appl.}\ 17:53--54, 1983
\newblock  (translated from the Russian paper in
  {\em Funktsional.\ Anal.\ i Prilozhen}).

\bibitem{E:tower} N.~D. Elkies.
\newblock  Explicit modular towers.
\newblock  In {\em Proceedings of the Thirty-Fifth
  Annual Allerton Conference on Communication, Control and Computing},
  pages 23--32.  Univ.\ of Illinois at Urbana-Champaign, 1998.
\newblock  {\sf http://arXiv.org/abs/math/0103107}

\bibitem{E:shimura} N.~D. Elkies.
\newblock  Shimura curve computations.
\newblock  In {\em Proceedings of ANTS-3}
  (Lecture Notes in Computer Science 1423), pages 1--47.
  Springer, Berlin, 1998.
\newblock  {\sf http://arXiv.org/abs/math/0005160}

\bibitem{E:drinfeld} N.~D. Elkies.
\newblock  Explicit towers of Drinfeld modular curves.
\newblock  In {\em Proceedings of the Third European Congress
  of Mathematics, Barcelona 2000}.
\newblock  {\sf http://arXiv.org/abs/math/0005140}

\bibitem{GS:drinfeld} A. Garcia and H. Stichtenoth.
\newblock  A tower of Artin-Schreier extensions of function fields
            attaining the Drinfeld-$\!$Vladut bound.
\newblock  {\em Invent.\ Math.}\ 121:211--233, 1995.

\bibitem{GS:more_drinfeld} A. Garcia and H. Stichtenoth.
\newblock  On the asymptotic behaviour of some towers of function fields
  over finite fields.
\newblock  {\em J.\ Number Theory} 61:248--273, 1996.

\bibitem{GS:classical} A. Garcia and H. Stichtenoth.
\newblock  Asymptotically good towers of function fields
    over finite fields.
\newblock  {\em C. R. Acad.\ Sci.\ Paris~I}\/ 322:1067--1070, 1996.

\bibitem{GST} A. Garcia, H. Stichtenoth, and M. Thomas.
\newblock  On towers and composita of towers of function fields
    over finite fields.
\newblock  {\em Finite Fields and their Appl.}\ 3:257--273, 1997.

\bibitem{Goppa} V.~D. Goppa.
\newblock  Codes on algebraic curves.
\newblock  {\em Soviet Math.\ Dokl.}\ 24:170--172, 1981.

\bibitem{GuSu} V. Guruswami and M. Sudan.
\newblock  Improved decoding of Reed-Solomon
    and algebraic-geometry codes.
\newblock  {\em IEEE Trans.\ Inform.\ Theory} 45:1757--1767, 1999.

\bibitem{Ihara} Y. Ihara:
\newblock  Some remarks on the number of rational points
    of algebraic curves over finite fields.
\newblock  {\em J. Fac.\ Sci.\ Tokyo} 28:721--724, 1981.

\bibitem{DiPippo} S.~A.~DiPippo.
\newblock  {\em Spaces of Rational Functions on Curves Over Finite Fields.}
\newblock  Ph.D. Thesis, Harvard, 1990.

\bibitem{Schanuel} S.~H. Schanuel.
\newblock  Heights in number fields.
\newblock  {\em Bull.\ Soc.\ Math.\ France} 107:433--449, 1979.

\bibitem{Serre} J.-P. Serre.
\newblock {\em Lectures on the Mordell-Weil Theorem} (trans.\ M.~Brown).
\newblock  F.~Vieweg \& Sohn, Braunschweig 1989.

\bibitem{ShWa} A.~M. Shokrollahi and H. Wasserman.
\newblock  List decoding of algebraic-geometric codes.
\newblock  {\em IEEE Trans.\ Inform.\ Theory} 45:432--437, 1999.

\bibitem{TV} M.~A. Tsfasman and S.~G. Vl\u{a}du\c{t}.
\newblock {\em Algebraic-Geometric Codes.}
\newblock Kluwer, Dordrecht, 1991.

\bibitem{TVZ} M.~A. Tsfasman, S.~G. Vl\u{a}du\c{t}, and T.Zink.
\newblock  Modular curves, Shimura curves and Goppa codes
    better than the Varshamov-Gilbert bound.
\newblock  {\em Math.\ Nachr.}\ 109:21--28, 1982.

\bibitem{Wan} D.~Wan.
\newblock  Heights and Zeta Functions in Function Fields.
\newblock  In {\em The Arithmetic of Function Fields}, pages 455-463.
  W. de Gruyter, Berlin, 1992.
   
\end{thebibliography}
\end{document}